\input amstex
\input epsf

 \documentstyle{amsppt}
 \magnification=1200
 \vsize=8.5truein
 \hsize= 6.0truein
 \hoffset=0.2truein
  \TagsOnRight

\nologo
 \NoRunningHeads
 \NoBlackBoxes

\topmatter
\title An optimal extension of Perelman's comparison theorem for quadrangles
and its applications
\endtitle

\author Jianguo Cao${}^1$, Bo Dai${}^2$ and Jiaqiang Mei${}^3$ \endauthor
\affil {University of Notre Dame, Peking University and Nanjing University} \endaffil
\address Department of Mathematics, University of Notre Dame,
Notre Dame, IN 46556, USA; Department of Mathematics, Nanjing University,
Nanjing 210093, China.
\endaddress
\email cao.7\@nd.edu \endemail
\thanks ${}^1$The first author is supported in part by an NSF grant. \endthanks
\address LMAM, School of Mathematical Sciences, Peking University, Beijing 100871, People's Republic of China
\endaddress
\email daibo\@math.pku.edu.cn \endemail
\thanks ${}^2$The second author is supported in part by NSFC Grant 10321001. \endthanks
\address Department of Mathematics and Institute of Mathematical Science,
Nanjing University,
Nanjing 210093, People's Republic of China
\endaddress
\email meijq\@nju.edu.cn \endemail
\thanks ${}^3$The third author is supported by Natural Science
Foundation of China (grant 10401015). \endthanks

%
\endtopmatter

\document
\baselineskip=.25in

\vskip1mm
\centerline{ Updated version, March 2009 }

\vskip2mm

\noindent {\bf Abstract:}  In this paper we discuss an extension of Perelman's comparison for quadrangles. Among applications of
 this new comparison theorem, we study the equidistance evolution of hypersurfaces in Alexandrov spaces with non-negative curvature.
 We show that, in certain cases, the equidistance evolution of hypersurfaces become totally convex relative to a bigger sub-domain.
 An optimal extension of 2nd variational formula for geodesics by Petrunin will be derived for the case of non-negative curvature.

     In addition, we also introduced the generalized second fundament forms for subsets in Alexandrov spaces. Using this new notion, we will
     propose an approach to study two open problems in Alexandrov geometry.



\vskip2mm

\head $\S 0$. Introduction
\endhead

In this paper, we derive a sharp comparison theorem for
quadrangles in complete metric spaces with non-negative curvature.
An earlier version of such a comparison for quadrangles was
discovered by Perelman with an asymptotical estimate (cf. \S 6 of
[Per91]). For smooth Riemannian manifolds with non-negative sectional
curvature, such sharp comparison theorem was implicitly stated in an important paper
[CG72] of Cheeger-Gromoll. We will extend Cheeger-Gromoll's approach to singular spaces of non-negative
curvature.

Among applications of our new sharp comparison for quadrangles, we
will derive a sharp version of the 2nd variational formula of
lengths in curved singular spaces, which improves the earlier work
of Petrunin (cf. [Petr98]). Other applications of our new sharp
comparison is to study the changes of Hessians of distance
functions in non-negatively curved spaces. This application would
provide a curved version of moving half-space method in Alexandrov
spaces with non-negative curvature, which we address in a separate
paper (cf. [CDM07]).

In \S 6 of an important preprint [Per91], Perelman pointed out
an asymptotic estimate for a class of quadrangles. Perelman's work
was completed before the notion of quasi-geodesic segments were introduced.
We will recall the definition of quasi-geodesics in \S 1 below. In what follows,
we always let $\sigma_{pq}:[0,\ell] \to M^n$ be a quasi-geodesic segment
from $p$ to $q$ of unit speed.

\proclaim{Definition 0.1}
Let $\triangle_1$ be a quasi-geodesic triangle with sides
$\{ \sigma_{qp} ,\sigma_{q\hat p}, \sigma_{p\hat p} \} $ and
$\triangle_2$ be a quasi-geodesic triangle with sides
$\{ \sigma_{q\hat p}, \sigma_{q\hat q}, \sigma_{\hat p \hat q} \} $.
If the two triangles have a common side $\sigma_{q\hat p}$ and if
$$
\angle_{q}(\sigma_{qp}'(0), \sigma_{q\hat p}'(0) ) + \angle_{q}
( \sigma_{q\hat p}'(0), \sigma_{q\hat q}'(0) )
= \angle_{q}(\sigma_{qp}'(0), \sigma_{q\hat q}'(0) )
\tag 0.1
 $$
 holds, then we say that the two quasi-geodesic triangles $\triangle_1$
 and $\triangle_2$ are co-planar at $q$.
 \endproclaim

 \proclaim{Theorem 0.2}(Extended Perelman's comparison for quadrangles)
  Let $M^n$ be a complete Alexandrov space
of $curv \ge 0$, $\triangle_1$ and $\triangle_2$ be two quasi-geodesic triangles  co-planar at $q$ as above. Suppose that
$\{ \sigma_{q\hat q},  \sigma_{p\hat p} \}$ are
the only two possible quasi-geodesics and the remaining quasi-geodesic segments $\{ \sigma_{\hat q \hat p},
\sigma_{q\hat p}, \sigma_{pq}
 \}$ are length-minimizing geodesic segments with
$$
\angle_{q}(\sigma_{qp}'(0), \sigma_{q\hat p}'(0) ) + \angle_{q}
( \sigma_{q\hat p}'(0), \sigma_{q\hat q}'(0) )
= \angle_{q}(\sigma_{qp}'(0), \sigma_{q\hat q}'(0) ) = \frac{\pi}{2}.
\tag 0.2
 $$
 Then

\smallskip
\noindent
(1) The function $s \to d(\sigma_{p\hat p}(s), \sigma_{q\hat q}(\Bbb R)) $ is concave at $s = 0$:
 $$
       \frac{d [f( \sigma_{p\hat p}(s)  )]}{ds}(0) \le \cos \theta, \tag 0.3
 $$
where $f(x) = d(x,\sigma_{q\hat q}(\Bbb R))$.

\smallskip
\noindent
(2)  If, in addition, $\sigma_{p\hat p}$ is a length-minimizing geodesic segment, then
 we have the following sharp estimates
 $$
d (  \hat p, \sigma_{q\hat q}(\Bbb R)) \le d(p, q) - s \cos \theta      \tag 0.4
 $$
where $s = d(p, \hat p)$ and
$\theta = \angle_{p}(\sigma_{p \hat p}'(0), \sigma_{pq}'(0) )$.
Moreover, the equality holds in (0.4) if and only if there is $\tilde q\in
\sigma_{q\hat q} (\Bbb{R})$ such that four points $\{ p,\hat p, \tilde q,
q\}$ span a totally geodesic flat trapezoid.
\endproclaim

For the special case of smooth Riemannian manifolds with non-negative sectional curvature,
the  comparison theorem for quadrangles above was due to Cheeger-Gromoll, see the proof of Theorem 1.10 of [CG72].

 Some independent work were carried out by Alexander-Bishop [AB08] via a different method. Among other things,
  Alexander and Bishop used the ratios of arc/chord and base-angle/chord to measure the convexity of
  hypersurfaces, which are very interesting. There are some overlap between our Corollary 2.4 and their results.
  On one hand, our hypothesis of corollary also allows points of zero geodesic
curvature for curves in surfaces, e.g., the super-graph of $y=x^4$ is strictly convex in
our sense but the boundary has zero geodesic curvature at the origin. Alexander-Bishop's definition of strictly convexity
 is not applicable to
case of  super-graph of $y=x^4$.
  On the other hand, Alexsnder-Bishop's work [AB08] covers not only spaces
with curvature $\ge k$ but also spaces with curvature $\le C$. Their results cover more ambient spaces.

\head $\S 1$. Earlier results on quasi-geodesics, development maps, classical comparisons and stability of Alexandrov spaces
\endhead

In this section, we review some earlier results for Alexandrov spaces, which will be used in later sections of our paper.
Other basic materials for Alexandrov spaces can also
be found in Chapter 10 of textbook [BBI01] and [BGP92].

A complete Alexandrov space $M^n$ is often referred to a complete metric space with curvature bounded below.  Since Alexandrov and his Russian
school of geometry used the geometric triangles to define ``the space with $curv \geq -k$", we require that $M^n$
must be a length space. Let $T_x^- M^n$ be the tangent cone of an Alexandrov space $M^n$ at $x$. With
some additional efforts, one can show that the tangent cone $T_x(M)$ must have angular measurement. In fact, the angular measurement is a
distance function defined on the unit tangent cone $\Sigma_x (M^n)$ of $M^n$ at $x$.
We now recall some important features of Alexandrov
spaces as follows.

\proclaim{Definition 1.1}   (1) A complete metric space $M^n$ is
called a length space (or an inner space) if for any pair of points
$\{ p,q\}$ in $M^n$, there exists a length-minimizing curve
$\sigma: [0,\ell]\to M^n$ from $p$ to $q$ such that the length
$L(\sigma)$ of $\sigma$ is equal to $d(p,q)$, the distance
between $p$ and $q$. Such a length minimizing curve $\sigma$
is called a geodesic segment of $M^n$.

\noindent (2) Suppose that $(M^n,d)$ is a complete length space.
We say that the space $(M^n,d)$ has $curv \geq 0$, if for any
length minimizing geodesic segment $\sigma: [0,\ell]\to M^n$
of unit speed and $p\not\in \sigma([0,\ell])$, the function
$$\eta_{p,\sigma}^0 (t) =\frac12 [d(p,\sigma(t))]^2
-\frac{t^2}{2} \tag 1.1
$$
is a concave function, i.e.
$$
\eta_{p,\sigma}^0 \left(\frac{t_1 +t_2}{2} \right)
\geq \frac12 [\eta_{p,\sigma}^0 (t_1) +\eta_{p,\sigma}^0 (t_2)],
\tag 1.2
$$
for all $\{ t_1,t_2 \} \subset [0,\ell]$.
\endproclaim

Roughly speaking, if the space $(M^n,d)$ has $curv \geq 0$,
then the triangle $\triangle_{p,\sigma}$ spanned by $\{ p,\sigma\}$
is more concave towards the vertex $p$ than the corresponding
triangle $\triangle_{p^\ast, \sigma^\ast }^\ast$ in the Euclidean
space $\Bbb{R}^2$.

Notice that
$$
[(t+c)^2]'' =(t^2)''
$$
for any constant $c$. Thus, in Definition 1.1 for the space with
$curv \geq 0$, we use the comparison $t^2$ instead of $(t+c)^2$
for appropriate $c$.

In order to define the Alexandrov space with $curv \geq -1$,
for any given pair $\{ p,\sigma\}$ with $p\not\in \sigma
([0,\ell])$, we choose the corresponding pair
$\{ p^\ast, \sigma^\ast \}$ in the hyperbolic plane
$\Bbb{H}^2$ of constance
curvature $-1$ more carefully as follows. Let
$$
\varphi_{p,\sigma} (t)=d(p,\sigma(t)).
\tag 1.3
$$
For any pair $\{ p^\ast,\sigma^\ast \} \subset \Bbb{H}^2$,
we let
$$
\varphi_{ p^\ast,\sigma^\ast}^\ast (t)=d_{\Bbb{H}^2}
(p^\ast,\sigma^\ast (t)),
$$
where $\sigma^\ast :[0,+\infty)\to \Bbb{H}^2$ is a geodesic of unit speed.
We require that
$$
\cases
\varphi_{p,\sigma} (0) =\varphi_{p^\ast,\sigma^\ast}^\ast (0), &  \\
\frac{d\varphi_{p,\sigma}}{dt} (0)
=\frac{d\varphi_{p^\ast,\sigma^\ast}^\ast}{dt} (0). &  \\
\endcases
\tag 1.4
$$

\proclaim{Definition 1.2} Let $(M^n,d)$ be a length space,
$\{ p,\sigma \}$, $\{ p^\ast,\sigma^\ast \} $,
$\varphi_{p,\sigma}$, $\varphi_{p^\ast,\sigma^\ast}$ be as above.
We say that the space $(M^n,d)$ has $curv \geq -1$ if for any
geodesic $\sigma: [0,\ell]\to M^n$ of unit speed and
$p\not\in \sigma ([0,\ell])$, the function $f_{-1}(t) = f_{p,\sigma,-1} (t)= - [1-\cosh
(d(p,\sigma(t))]$ satisfies  the differential inequality:
$$
 f_{-1}'' \le [1 + f_{-1}].   \tag 1.5
$$

Similarly, for $k=1$, we let $\rho_1(s)=[1-\cos s]$,
$f_{p,\sigma,1} (t) =\rho_1(d(p,\sigma(t))$.
If $f_1(t) = f_{p,\sigma,1} (t)$ satisfies the differential inequality:
$$
 f_{1}'' \le [1 - f_{1}] \tag 1.6
$$
for all length-minimizing geodesic $\sigma$, then we say that $M^n $ has $curv  \ge 1$.
\endproclaim

It is well-known that if $(M^n,d)$ has $curv \geq k$, then
$(M^n,\lambda d)$ has $curv \geq \frac{k}{\lambda^2}$, where
$\lambda >0$. By scaling the distance function with a factor
$\lambda >0$, we can define the notion ``$curv \geq k$" for any
real number $k\in \Bbb{R}$.

Burago, Gromov and Perelman [BGP92] derived several important
results for Alexandrov spaces. Among other things, they discovered
the following results.

\proclaim{Theorem 1.3}([BGP92, \S 7.8.1])
Let $(M^n,d)$ be an Alexandrov space with $curv \geq k$. Then

\noindent (1.3.1) The dimension of $(M^n,d)$ is an integer or infinity,
which is equal to the Hausdorff dimension;

\noindent (1.3.2) Let $T_p^- (M^n)$ be the cone over the space
of directions of $M^n$ at $p$, and let $sM^n =(M^n,s\,d)$ be the
scaling of the space  $(M^n,d)$. Then $T_p^- (M^n)$ is isometric to
the Gromov-Hausdorff limit of the pointed spaces
$\{ (sM^n,p)\}$ as $s\to +\infty$;

\noindent (1.3.3) The tangent cone $T_p^- (M^n)$ has $curv \geq 0$;

\noindent (1.3.4) Let $\Sigma_p(M^n) =\{ \vec{v}\in
T_p^- (M^n) : |\vec{v}|=1 \}$. Then the unit tangent cone
$\Sigma_p(M^n)$ has $curv \geq 1$.
\endproclaim

There are many non-smooth spaces with $curv \geq  0$.

\bigskip
\noindent
{\bf Example 1.4} (i) Let $M^n$ be a complete
smooth Riemannian manifold with $curv \geq k$.
Suppose $\Omega\subset M^n$ is a convex sub-domain
of $M^n$, where $\Omega$ is not necessarily smooth.
Then, by a theorem of Buyalo, the boundary $\partial
\Omega$ of $\Omega$ has $curv \geq k$ as well.

\noindent (ii) $M^2=\{ (x,y,z)\in \Bbb{R}^3|
z=\sqrt{x^2+y^2} \}$ has $curv \geq 0$, but $(0,0,0)=p$
is not a smooth point of $M^2$.

\noindent (iii) Let $\Omega$ be an American football.
Its boundary $M^2=\partial\Omega$ has $curv \geq 0$.

\noindent (iv) If $\Sigma$ is an Alexandrov space with
$curv \geq 1$ (e.g. $\Sigma  =\Bbb{C} P^m$ is a complex
projective space with the Fubini-Study metric), then we consider
the cone over $\Sigma$ as follows. Let $M=Cone_0 (\Sigma)
=\Sigma\times [0,+\infty)/\Sigma\times\{ 0\}$.  For any pair of
points $\{ (p,t_1),(q,t_2)\} \subset \Sigma\times [0,+\infty)$,
we define
$$
[d_M ((p,t_1),(q,t_2))]^2 =t_1^2 +t_2^2 -2t_1 t_2 \cos [d_\Sigma
(p,q)]. \tag 1.6
$$
It was shown in [BGP92] that $(M,d_M)$ above has $curv \geq 0$.

\bigskip

 Recall that if for any
length-minimizing geodesic $\sigma: [0, \ell] \to M^n$ of unit speed we have
$$
\frac{d^2 f(\sigma(t))}{dt^2} (0) \le c
$$
in barrier sense then we say
$$
Hess(f) (\sigma'(0), \sigma'(0) ) \le c
$$
in barrier sense. With some additional efforts, one can easily verify the following estimate for the
upper bound of Hessian of distance functions.

\proclaim{Proposition 1.5} Let $(M^n,d)$ be a complete Alexandrov space
with $curv \geq k$. Suppose that  $\varphi: [0, d] \to M^n$ is a length-minimizing geodesic of unit speed
from $x$ to $p$, and that $\sigma: [0, \epsilon] \to M^n $ is a quasi-geodesic of unit speed
with $\angle_x(\sigma'(0), \varphi'(0)) = \frac{\pi}{2}$.  Then, in barrier sense, for $d_p(x)=d(p,x)$ we have
$$
\roman{Hess} (d_p) (\sigma'(0), \sigma'(0) ) \leq
\cases \frac{1}{d_p}, & \roman{if} \ k=0, \\
 \cot \left( {d_p} \right) ,
& \roman{if}\ k =1,\\
 \coth \left( {d_p}  \right).
& \roman{if}\ k = -1.\\
\endcases
\tag 1.7
$$
\endproclaim
\demo{Proof} It is clear that
$$
    \frac{d^2[ h(u(t))]}{dt^2 }  =  h'(u(t)) u''(t) + h''(u(t)) [u'(t)]^2. \tag *
$$

Let $u(t) = d(p, \sigma(t))$.
   Thus, by our assumption we have  $u'(0) = 0$.

 (1)   When $k = 0$, we let $h(u) = \frac{u^2}{2}$.

    By definition, we see that if $f(t) = h(u(t))$ then $f''(t) \le 1$. It follows that
   $$
  1 \ge f''(t) = u(t) \roman{Hess} (d_p) (\sigma'(0), \sigma'(0) ) + 0
   = d_p(t) \roman{Hess} (d_p) (\sigma'(0), \sigma'(0) ).
   $$
It follows that $\roman{Hess} (d_p) (\sigma'(0), \sigma'(0) ) \le \frac{1}{d_p}$ for $k = 0$.

  \bigskip

  (2)  When $k \ge 1$, we let $h(u) = 1 - \cos u$, $u(t) = d(p, \sigma(t))$ and $f(t) = h(u(t))$.
It is known that $ f''(t) \le [ 1 - f(t)]$. It follows that
$$
   - \frac{d^2 [ \cos d_p(t)]   }{dt^2} \le \cos d_p(t). \tag \ddag
$$

   Therefore, by (*)-($\ddag$) and initial condition $u'(0) = 0$,  we have
   $$
[\sin d_p(t)] \roman{Hess} (d_p) (\sigma'(0), \sigma'(0) ) \le \cos d_p(t).
   $$

\bigskip

   (3) The case of $k = -1$ can be handled similarly. \qed
\enddemo

When $(M^n,d)$ has $curv \geq k >0$, then the diameter is less
than or equal to $\frac{\pi}{\sqrt{k}}$, i.e.
$$
\roman{Diam} (M^n) \leq \frac{\pi}{\sqrt{k}}.
\tag 1.8
$$
Thus, the estimate (1.7) above makes sense for all $k>0$ as well.

    There is another way to see why the Hessian inequality (1.7) above holds in barrier sense, by the development maps used by Russian school of
    geometry, see \S 7.3 of Plaut's survey paper [Pl02, p861]. In fact, there are several  equivalent definitions of  quasi-geodesics.

    Let us first recall a simple definition of quasi-geodesics without using development map.

\proclaim{Definition  1.6} (Quasi-geodesics, [Pl02, p860]) Let $Y$ be a metric space,   $\eta : [a, b] \to Y$ be a Lipschitz curve of
unit speed and let $d_q(t) = d(q, \eta(t))$. Then $\eta$ is called a quasi-geodesic if for every $q \in Y$ there is a function $h$ with
$\lim_{s \to 0^+}
\frac{h(s)}{s} = 0$, and $[d_q(t)^2 - t^2]'' \le h(d_q(t))$ for all $t \in [a, b]$. We write $[d_q(t)^2 - t^2]'' \le o(d_q(t))$
for short.
\endproclaim

    Using the triangle comparison theorems for Alexandrov space with $curv \ge k$, one can show that

\proclaim{Proposition 1.7} Let $q \in Y$, $\eta : [a, b] \to Y$ be a geodesic of unit speed and let $d_q(t) = d(q, \eta(t))$.
Suppose that $Y$ has $curv \ge k$ and that
$$
f(t) =
\cases \frac 12 d^2_q(t), & \roman{if} \ k=0, \\
 \frac{1- \cos (\sqrt{k} d_q(t)) } {k} ,
& \roman{if}\ k>0,\\
\frac{1- \cosh (\sqrt{-k} d_q(t)) } {k}.
& \roman{if}\ k<0.\\
\endcases
$$
Then
$$
      f''(t) \le [ 1 - k f(t)]. \tag 1.9
$$

Moreover, the above condition holds if and only if the conclusion of Toponogov comparison theorem holds for
any geodesic hinges in $Y$.

\endproclaim

Perelman and Petrunin (cf. [PP94] and [PP96]) used the inequality (1.9) to define quasi-geodesics,
see Definition 2.1 below as well.

The third definition of quasi-geodesics uses the development maps for a curve in the model space $M^2_k$,
where $M^2_k$ is a complete
    simply-connected surface
    of constant curvature $k$. We thank Professor Stephanie Alexander for supplying us a correct definition of
    development maps.

\proclaim{Definition 1.8} ([AB96] Development map of a curve relative to
p) Let $\sigma_{xy}$ be a length-minimizing curve of unit speed in a metric space.
Suppose $\gamma$ is a rectifiable curve parameterized by arc-length
in $M$ and $p$ is a point not on $\gamma$.  The {\it Alexandrov
development} $\tilde\gamma$ of $\gamma$ from $p$ is a curve in $M^2_k$
obtained as follows: associate to $p$ a point $\tilde p$ in $M^2_k$, and
associate to the minimizer $\sigma_{p\gamma(t)}$ to a minimizer
$\tilde\sigma_{\tilde p \tilde \gamma(t)}$ of the same length, turning monotonically in $t$ in
such a way that $t$ is also the arc-length parameter of $\tilde \gamma$. The union of the
$\tilde \sigma_{\tilde p \tilde \gamma(t)}
$ is called the
{\it ``cone of the development"}. A development is possible whenever
the maximum distance from $p$ to $\gamma$ is less than $\pi/
\sqrt{  \max\{0, k\}}$.

\endproclaim

We now state another equivalent definition of quasi-geodesics.

\proclaim{Proposition 1.9} (page 861 of [Pl02] characterization of quasi-geodesics via development maps). Let $Y$ be a
complete Alexandrov space with
$curv \ge k$ and let $c: [a, b] \to Y$ be a curve of unit speed. Suppose that $\tilde c^*$ is a $k$-development
    of $c$ at $p$ with $p  \notin c([a, b])$.
Then $c$ is a quasi-geodesic if and only if $c^*$ is $k $-convex in following sense:

\smallskip
\noindent
(1.9.1) (Alexandrov) For any $p \in M$ and geodesic $\gamma: [a, b] \to M$, the $k$-development
of $\gamma$ at $p$ is convex, where $p \notin \sigma([a, b])$ and $M$ has curvature $\ge k$.

\smallskip
\noindent
(1.9.2) Suppose that $M$ has curvature $\ge k$. A curve $c$ in $M$ is a quasi-geodesic if and only
if it has unit speed and its development of $\tilde c$ is $k$-convex relative to any point $p \notin c$.
\endproclaim

It should be pointed out, for some non-smooth Alexandrov space $(M^n,d)$,
the distance function $d_p(x) =d(p,x)$ may {\it not\/} be convex for
$x\in B_\varepsilon (p)$, where
$B_\varepsilon (p)=\{ q| d(p,q)<\varepsilon \}$.
More precisely,  the conclusion of the following Proposition
for smooth Riemannian manifolds might fail for Alexandrov spaces.

\proclaim{Proposition 1.10} Let $(M^n,d)$ be a complete
$C^2$-smooth Riemannian manifold with $curv \geq k$.
Suppose that $\overline{\Omega}$ is a compact subset of $M^n$.
Then there exists $\varepsilon_0 >0$ depending only on
$\overline{\Omega}$ such that $B_\varepsilon (p)$ is a convex subset
for all $p\in \overline{\Omega}$ and $0<\varepsilon \leq
\varepsilon_0$.
\endproclaim

Here is a simple example of non-smooth Alexandrov space $M^2$
for which the conclusion of Proposition 1.10 fails.

\proclaim{Example 1.11} Let $M^2=\{ (x,y,z)\in\Bbb{R}^3 |
z=\sqrt{x^2+y^2} \}$ and $\overline{\Omega} =\{
(x,y,z)\in M^2| 0\leq z\leq 1 \}$. We choose a sequence of
points $p_\varepsilon =(\varepsilon,0,\varepsilon)$. It is clear
that the metric disk $B_{2\varepsilon} (p_\varepsilon)$ is not
convex for $0 <  \varepsilon\leq 1$.
\endproclaim

The above example indicates that $d_{p_0}: B_\varepsilon (p_0)
\to \Bbb{R}$ is not necessarily a convex function if there is a sequence
of singular points $\{ q_i\} \to p_0$ as $i\to +\infty$. In addition,
the injectivity radius of $M^n$ restricted to  $\overline{\Omega}$
is zero since
$$\roman{inj}_{M^n} (p_\varepsilon)\leq 2\varepsilon \to 0
$$
as $\varepsilon\to 0$.

In order to show that the distance function $d_p(x)=d(p,x)$ is
free of critical points on $[B_\varepsilon (p) -\{ p\}]$ in the
sense of Gromov, Perelman cleverly constructed a local convex function
$\psi_p: B_\delta (p) \to \Bbb{R}$ which is bi-Lipschitz comparable to
$d_p$ as follows.

\proclaim{Theorem 1.12}(Perelman [Per94b], Kapovitch [Ka02, Lemma 4.2])
Let $M^n$ be a finite dimensional Alexandrov space with
$curv \geq k$. For each $p\in M^n$, there is a $\delta =
\delta(p)$ depending only on the local volume growth
of $B_\varepsilon (p)$ and there is a strictly convex non-negative
function $\psi_p: B_\delta (p)\to [0,+\infty)$
with $\psi_p(p)=0$ and
$$
B_{\frac{\varepsilon}{\lambda}} (p) \subset
\psi^{-1}_p ([0,\varepsilon])
\subset B_{\lambda\varepsilon} (p)
$$
for some $\lambda\geq 1$ and $\varepsilon\in
(0,\frac\delta\lambda ]$.
Consequently, the distance function $d_p$ has no
critical point on $[B_{\frac{\delta}{\lambda}}(p)
-\{ p\}]$ in the sense of Gromov.
\endproclaim

Using Theorem 1.12 above, Perelman [Per94] further studied
the local structure of Alexandrov spaces. In Example 1.4 above,
we see that the cone over an Alexandrov space $\Sigma^{n-1}$ with
$curv \geq 1$ has $curv \geq 0$. In [BGP92], Burago-Gromov-Perelman
also constructed parabolic or hyperbolic cones over a lower
dimensional Alexandrov space $\Sigma^m$, the resulting cones have
$curv \geq k_1$ for some other $k_1 \in \Bbb{R}$.

\proclaim{Definition 1.13}(Perelman's MCS-spaces)
We define MCS-spaces (spaces with multiple conic singularities)
inductively on the dimensions.

(0)  A point $\{ p\}$ is a $0$-dimensional MCS-space;

(1) We say that $M^n$ is an MCS-space if for each $p\in M^n$,
there is a small ball $B_\varepsilon (p)$ of radius $\varepsilon$
centered at $p$ such that $B_\varepsilon (p)$ is homeomorphic
to a cone over a lower dimensional MCS-space $\Sigma$.
\endproclaim

A remarkable result of Perelman asserts that all possible singularities of
any finite dimensional Alexandrov space $M^n$ are at most as bad as conic
singularities.

\proclaim{Theorem 1.14}(Perelman [Per94b])
Let $M^n$ be a complete finite dimensional Alexandrov space with
$curv \geq k$. Then $M^n$ must be an MCS-space.
\endproclaim

The proof of Perelman's structure theorem above is inspired by the
celebrated Perelman's stability theorem.

\proclaim{Theorem 1.15}(Perelman [Per94b], [Ka07], [BBI01,p400]) For any $k\in \Bbb{R}$ and any compact Alexandrov space $M^n$ with $curv \geq
k$, there is an $\varepsilon =\varepsilon(M^n)>0$ such that every compact Alexandrov space $Y^n$ with $curv \geq k$,
$$
d_{GH} (M^n,Y^n)<\varepsilon
$$
and $\dim (M^n)=\dim (Y^n)=n$ must be homeomorphic
to $M^n$, where $d_{GH}(X,Y)$ denotes the Gromov-Hausdorff
distance between $X$ and $Y$.

Consequently, $B_\varepsilon (p)$ is homeomorphic to a small
ball $B_\varepsilon (O_p)$ in the tangent cone $T_p^-
(M^n)$ for sufficiently small $\varepsilon >0$.
\endproclaim

For a non-compact complete space $M^n$, one consider
$\{ (M_i^n,p_i)\} \to (M^n,p)$ in the  pointed Gromov-Hausdorff
convergence. A Similar version of Theorem 1.15 for a sequence of
pointed Alexandrov spaces also exists, see [Ka07].

Perelman's stability theorem can be used to simplify several crucial steps
in Perelman's solution to Thurston's geometrization conjecture.

\head $\S 2$. Proof of the new sharp comparison theorem for
quadrangles
\endhead

In this section, we will provide a detailed proof of Theorem 0.2.

Let us first recall an equivalent definition of quasi-geodesics due to
Perelman and Petrunin.  Let
$$
\rho_k(x)=
\cases \frac{x^2}{2}, & \roman{if} \ k=0; \\
\frac1{k} [1-\cos ( \sqrt{k}x ) ],
& \roman{if}\ k>0;\\
\frac1{k} [1-\cosh ( \sqrt{-k}x )] ,
& \roman{if}\ k<0.\\
\endcases
$$
One considers
$$
f_{p,\sigma,k} (t) =\rho_k(d(p,\sigma(t)).
\tag 2.1
$$

\proclaim{Definition 2.1} ([PP94], [PP96]) Let $(M^n,d)$ be a complete
Alexandrov space with $curv \geq k$, and let $\sigma:[0,\ell] \to M^n$ be a
Lipschitz curve of unit speed. If
$$
f''_{p,\sigma,k}(t)\leq 1-k f_{p,\sigma,k}(t)
$$
holds for any $p\in M^n$, then $\sigma$ is called
a quasi-geodesic segment.
\endproclaim

We remark that one might be able to use the generalized
2nd fundamental form for curves to provide an equivalent
definition of quasi-geodesics, see \S 4 below.

It is known that any geodesic segment is a quasi-geodesic.
Petrunin also observed that comparison theorem holds for a
class of quasi-geodesic hinges.

\proclaim{Proposition 2.2} Let $\sigma_1:[0,\ell_1]\to M^n$
be a length-minimizing geodesic and
$\sigma_2:[0,\ell_2]\to M^n$ be a quasi-geodesic with
$p=\sigma_1(0)=\sigma_2(0)$ and
$$\theta =\angle_p (\sigma'_1 (0), \sigma'_2 (0) ).
$$

\smallskip
\noindent
(1) Suppose that $M^n$ has $curv \geq 0$. Then
$$
[d(\sigma_1(\ell_1),\sigma_2(\ell_2))]^2 \leq
\ell_1^2 +\ell_2^2  -2\ell_1 \ell_2 \cos\theta.
\tag 2.2
$$

\smallskip
\noindent
(2) Suppose that $M^n$ has $curv \geq 1$. Then
$$
\cos [d(\sigma_1(\ell_1),\sigma_2(\ell_2))] \ge
(\cos \ell_1) (\cos \ell_2) + (\sin \ell_1)   (\sin \ell_2) \cos \theta.
$$
\endproclaim
\demo{Proof} (1) When $k = 0$, we let $f(t) = \frac{ [d( \sigma_2 (t), \sigma_1(\ell_1)) ]^2    }{2} $.
By an equivalent definition of quasi-geodesics, we have $f''(t) \le 1$. Using the initial condition
$f(0) = \frac{\ell_1^2}{2} $ and $f'(0) = \ell_1 \cos \theta$, we see that $f''(t) \le 1$ implies
$$
   2 f(t) \le \ell_1^2 + t^2  - 2 t \ell_1 \cos\theta.
$$
The inequality (2.2) follows.

(2) When $k = 1$, we will use an observation of Gromov to cancel the first derivatives. Let
$h(t) =  \cos [d(\sigma_2 (t), \sigma_1(\ell_1))]$ and $h^*(t) =
(\cos t ) (\cos \ell_1) + (\sin t )   (\sin \ell_1) \cos \theta$.

    When $f(t) = 1 -h(t)$, by an equivalent definition of quasi-geodesics, we have $f''(t) \le [ 1 - f(t)]$.
    It follows that
    $$  h''(t) \ge h(t). $$

    Inspired by Gromov, we let
    $$\eta(t) = h'(t) h^*(t) - [h^*(t)]' h(t).$$

    By the inequality $h''(t) \ge h(t)$, we see that
    $$
\eta'(t) \ge 0,
    $$
whenever $\min\{h(t), h^*(t)\} \ge 0$.

By our assumption, we see that $\eta(0) = 0$.  It follows from $\eta'(0) \ge 0$ that
$\eta(t) \ge 0$.  Whenever $\min\{h(t), h^*(t)\} \ge 0$, we have
$
\{\log [h(t)]\}' \ge \{\log [h^*(t)]\}'
$
and hence
$$
  h(t) \ge h^*(t).
$$
The other cases could be similarly, we leave it to readers.  In fact,  the differential inequality $f''(t) \le [1 -f(t)]$ with the initial
conditions above would implies $h(t) \ge h^*(t)$, see textbook [Pete98], page 327-330. \qed
\enddemo

\bigskip
\bigskip
\noindent
\demo{\bf Proof of Theorem 0.2}

We first verify Theorem 0.2 (2).

Let $\alpha =\angle_q (\sigma'_{q\hat p} (0),
\sigma'_{q\hat q} (0))$,  $\tilde q=
\sigma_{q\hat q}  (L\cos\alpha )$ and $L=d(q,\hat p)$.
Since $M^n$ has $curv \geq 0$ and $\sigma_{q\hat p}$
is length-minimizing, we have
$$
[d(\hat p,\tilde q)]^2 \leq L^2 +(L\cos\alpha)^2 -2(L\cos \alpha)^2
=(L\sin\alpha)^2.
\tag 2.3
$$
Let $T=d(p,q)$. Applying comparison theorem twice to the
geodesic triangle $\triangle_1$, we also obtain
$$
L^2 \leq [T^2 +s^2-2sT \cos\theta] \tag 2.4
$$
and
$$
\align
s^2 &\leq \left[L^2 +T^2-2LT\cos \left(\frac{\pi}{2}
-\alpha\right) \right] \\
& =L^2 +T^2 -2LT \sin\alpha. \tag 2.5
\endalign
$$
Using (2.3)-(2.5), we have
$$
\align
d(\hat p,\tilde q) &\leq L\sin\alpha \\
&\leq \frac{1}{2T} [L^2 +T^2 -s^2] \\
&\leq \frac{1}{2T} [(T^2 +s^2-2sT\cos\theta) +(T^2 -s^2)] \\
&\leq \frac{1}{2T} [2T^2 -2sT\cos\theta] \\
&=T-s\cos\theta.
\endalign
$$
It follows that
$$
\align
d(\hat p, \sigma_{q\hat q} (\Bbb{R})) &\leq d(\hat p,
\tilde q) \\
& \leq T-s\cos\theta.
\endalign
$$
This completes the proof of Theorem 0.2 (2).

To see that the first assertion of Theorem 1.2 is true, we use the development maps of $\sigma_{p\hat p}$  and
$\sigma_{q\hat q} $ relative to the midpoint
$q_{mid} = \sigma_{pq}(\frac T2)$.

 We see that, by (1.9.2) that $\sigma^*_{p, \hat p}$ lies {\it inside } of the trapezoid in the
model space $M^2_0 = \Bbb R^2$. We now use the fact $\sigma^*_{p, \hat p}$ has unit speed to conclude that
$$
       d(\sigma^*_{p, \hat p}(t), \sigma^*_{q\hat q} (\Bbb{R})) \le  d(\varphi^*(t'), \sigma^*_{q\hat q} (\Bbb{R})) \tag2.6
$$
with $\frac{t'}{t} \to 1$ as $t \to 1$, where $\varphi^*: [0, \infty) \to \Bbb R^2$ is a straight line with $\varphi^*(0) = p$
and $[\varphi^*]'(0) = [\sigma^*_{p, \hat p}]'(0)$.

   In the model space, we have
$$
d(\varphi^*(t'), \sigma^*_{q\hat q} (\Bbb{R})) \le T - t' \cos \theta. \tag2.7
$$
It follows from (2.6)-(2.7) and discussion above that the inequality (0.3) holds in barrier sense.

The equality case in (0.4) holds if and only if the four points span a totally geodesic flat trapezoid in $M^n$. \qed
\enddemo

In his important preprint [Per91], Perelman showed that, for fixed
$\theta$, one has
$$
d(\hat p, \sigma_{q\hat q} (\Bbb{R}))  \leq T-s\cos\theta
+o(s^2),
$$
where $\lim\limits_{s\to 0} \frac{o(s^2)}{s^2} =0$.

We present some direct applications of Theorem 0.2 and its proof to the equidistance evolutions.

\proclaim{Definition 2.3}  Let $\Omega$ be a subset of a complete Alexandrov space $M^n$ with $curv \ge k$.

(1) If for any pair of points $\{p, q\} \subset \Omega$ there is a length-minimizing geodesic $\sigma_{p, q}: [0, \ell]
\to M^n$ from $p$ to $q$ such that $\sigma((0, \ell)) \subset \Omega$, then $\Omega$ is called convex.

(2) Let $\Omega$ be a convex domain and let $\overline{\Omega}$ be its closure.
If any quasi-geodesic $\sigma: [0, \ell] \to M$ tangent to the boundary $\partial \Omega$ of
$\Omega$ at $\sigma(0)$
 has the property that $\sigma(t) \notin \overline{\Omega} $ for $t \in (0, \delta]$ and some positive $\delta \le \ell$,
 then $\Omega$ is called a {\bf strictly} convex domain.

(3) Let $\Omega \subset U$. If for any pair of points $\{p, q\} \subset \Omega$ and any length-minimizing geodesic $\sigma_{p, q}: [0, \ell]
\to U$ from $p$ to $q$ such that $\sigma((0, \ell)) \subset \Omega$, then $\Omega$ is called totally convex relative to $U$.

\endproclaim

\proclaim{Corollary 2.4}
Let $M^n$ be a complete Alexandrov space with $curv \geq 0$.
Suppose that $\Omega_0$ is a compact convex domain in
$M^n$ and $\partial \Omega_0$ is strictly convex. Then
$\Omega_{-T} =\{ p\in\Omega_0 | d(p,\partial \Omega_0)
\geq T \}$ is strictly convex
for all $T<T_0 =\max \{ d(p,\partial\Omega_0) | p\in
\Omega_0 \}$.
Furthermore, $\Omega_{-T}$ is totally convex in $\Omega_0$,
$\dim [\Omega_{-T_0}] =0$ and $\Omega_0$ is
contractible.
\endproclaim

  \medskip

  \noindent
  {\bf Remark:} For the special case when $M^n$ is a complete smooth Riemannian manifold with non-negative curvature but positive curvature outside
  a compact set, Corollary 2.4 was implicitly proved by Cheeger-Gromoll, see page 421 and page 431 of [CG72].

   \smallskip

\demo{Proof of Corollary 2.4} We present a proof inspired by Cheeger-Gromoll [CG72] with some modifications and will use the proof of Theorem 0.2.

  Let $f: M^n \to \Bbb{R} $ be a signed distance as follows: $ f(x) = d(x, \partial \Omega) $ for $x \in  \Omega$ and
  $f(x) = - d(x, \partial \Omega) $ for $ x \notin \Omega$.

\medskip
\noindent
{\bf Step 1. Proof of weak convexity.}
\smallskip

  We first consider the signed distance function $f(x)$ restricted to a length-minimizing geodesic segment in the same way as Cheeger-Gromoll did in [CG72].
  Let $\sigma: [0, \ell] \to \Omega$ be a length-minimizing geodesic of unit speed. For each interior point $p$ of $\Omega$, we let
$\Lambda_{p, \partial\Omega} = \{ \vec v \in T^-_p(M) \quad |
\quad \vec v = \varphi'(0), \varphi: [0, f(p)] \to M, |\varphi'(t)
| = 1, \varphi(0) = p, \varphi(f(p)) \in \partial \Omega\}$ be the
subset of all unit minimizing directions from $p$ to $\Omega$.
Since $Min_{p, \partial\Omega}$ is a compact subset of $\Sigma_p$,
there are $q \in \partial \Omega$ and $\vec w = \sigma_{pq}'(0)
\in \Lambda_{p, \partial\Omega}$ such that
$$
\theta = \angle_p(\sigma'(0), \vec w ) = \inf\{
\angle_p(\sigma'(0), \vec v)    \quad  | \quad  \vec v \in Min_{p,
\partial\Omega} \}
$$
and
$$
   d(p, q) = f(p) = d(p, \partial \Omega)
$$
hold.

   Let us verify the following.

\smallskip
 {\bf Claim 2.5.} Let $p \in int(\Omega)$, $q \in \partial \Omega$, $f: M \to \Bbb R$, $\theta$ and $\sigma:
 [0,\ell] \to \Omega$ be as above. Then, for each $s> 0$, there is a quasi-geodesic $\psi: [0, \delta] \to M$ tangent to
 $\partial \Omega$ at $\psi(0) = q$ such that
 $$
  f(\sigma(s)) = d( \sigma(s), \partial \Omega) < d( \sigma(s), \psi(t_s)) \le f(p) - s \cos \theta \tag2.8
 $$
 holds for some $t_s > 0$.
\smallskip

Recall that $d(p, q) = d(p, \partial \Omega)$, by the first variational formula, our length-minimizing geodesic
$\sigma_{pq}$ is orthogonal to $\partial \Omega$ at $q$. As Perelman observed that the tangent cone
$T^-_q(\overline{\Omega})$ has metric splitting
$$
T^-_q(\overline{\Omega})   = [0, \infty) \times T^-_q(\partial \Omega). \tag2.9
$$

 We now choose a length-minimizing geodesic segment from $q$ to $\sigma(s)$ of unit speed, say $\sigma_{q\sigma(s)}: [0, L]
 \to \Omega$. Let $\vec \xi  = \sigma_{q\sigma(s)}'(0)$. Using (2.9), we can write
 $$
\vec \xi = (\cos \alpha) \vec \eta + (\sin \alpha) \vec w^*  \tag2.10
 $$
 for some $ \vec \eta \in T^-_q(\partial \Omega)$  and $ 0 < \alpha < \frac{\pi}{2}$, where $\vec w^*$ is the left derivative of
 $\sigma_{pq}: [0, d(p, q)] \to M$ at its endpoint $q$ and $| \vec \eta| = 1 = |\vec w^*|$.

  For any unit direction $\vec \eta$, there is a quasi-geodesic $\psi: [0, \delta] \to M$ with
  $\psi'(0) = \vec \eta$. Let  $t_s = L \cos \alpha$ where $L = d(q, \sigma(s))$. By the proof of Theorem 0.2, we
  have
  $$
    d(\sigma(s), \psi(t_s)) \le d(p, q) - s \cos \theta. \tag2.11
  $$
Because $t_s = L \cos \alpha > 0$ and our domain $\Omega$ is strictly convex, by definition 2.3 (2), we have
$$
  f(\sigma(s)) =  d(\sigma(s), \partial \Omega) < d(\sigma(s), \psi(t_s)). \tag2.12
$$

   Claim 2.5 follows from (2.11)-(2.12).

    It now follows from Claim 2.5 that $f(x) = d(x, \partial \Omega)$ is a concave function for $x \in \Omega$.

\medskip
\noindent
{\bf Step 2. Proof of strict convexity.}
\smallskip
    Let
    $\Omega_{-c} = f^{-1}([c, \infty))$ for $c\ge 0$.  We would like to show that
    the convex domain $\Omega_{-c}$ is strictly convex for $c >0$ by using (2.8) and a theorem of
    Perelman-Petrunin. Let $\phi: [0, \delta^*] \to M$ be a quasi-geodesic
    segment tangent to $\partial \Omega_{-c}$ at $p$. Our goal is to verify
    $$
        f( \phi(s)) < f(\phi(0)) \tag**
    $$
for any $s \in (0, \delta]$ and some $\delta >0$.

   For special case when the above quasi-geodesic
    $\phi: [0, \delta^*] \to M$ is a length-minimizing geodesic segment, the inequality
$f( \phi(s)) < f(\phi(0))$ is a direct consequence of (2.8) by choosing $\theta = \frac{\pi}{2}$.

    For the general case of tangential quasi-geodesic $\phi: [0, \delta^*] \to M$, we use Corollary 3.3.3 of [Petr07] to get
    the fact that the function
    $$
    \eta(s) = f(\phi(s)) \tag2.13
    $$
is a concave function of $s$ with $\eta'(0) = 0$. It follows from concavity and $\eta'(0) = 0$ that
$$
\eta(s) = f(\phi(s)) \le f(\phi(0)) \tag2.14
$$
for all $s \ge 0$.

   If $f(\phi(s_0)) <  f(\phi(0))$ for some $s_0 >0 $, by concavity we have that
   $$
     f(\phi(s)) \le  f(\phi(s_0)) + (s-s_0) \frac{f(\phi(s_0)) -  f(\phi(0))}{s_0} \le f(s_0), \tag2.15
   $$
for all $s \ge s_0$.

Hence, we may assume that there were $s_0 >0 $ such that
$$
f(\phi(s_0)) = f(\phi(0))  \tag2.16
$$
and we will derive a contradiction as follows.

   Choose a length-minimizing geodesic $\sigma_{pq}: [0, c] \to M$ from $p$ to $\partial \Omega$ such that
   $q \in \partial \Omega$, $d(p, q) = c = d(p, \partial \Omega)$ and
   $$
 \angle_p(\phi'(0),  \sigma_{pq}'(0)   ) = \inf\{\angle_p(\phi'(0), \vec w  ) \quad | \quad \vec w
\in Min(p, \partial\Omega)  \} = \frac{\pi}{2}.
   $$

   We further choose the midpoint $q_{mid} = \sigma_{pq}(\frac c2)$. We will consider the development map of the quasi-geodesic
   $\phi$ relative to $q_{mid}$. Since $q_{mid}$ is an interior point of a length-minimizing geodesic $\sigma_{pq}$,
   its log image $log_p( q_{mid} ) $ is unique, which is equal to $c\sigma_{pq}'(0) = \tilde q_{mid}$. Let us consider a ``half-space"
   in the tangent cone $T_p(M)$ relative to $\tilde q_{mid}$ as follows:
$$
    Half_{p, \tilde q_{mid}} = \{ \vec u \in T^-_p(M) \quad | \quad  \langle \vec u, \tilde q_{mid}) \ge 0 \}.
$$

     Because $\phi: [0, s_0] \to M$ is a quasi-geodesic, its ``cone of development"  relative $q_{mid}$ can be isometrically embedded
     into $Half_{p, \tilde q_{mid}}$ with vertex $\tilde q_{mid}$. It follows that there is $\vec u_0 \in log_p(\phi(s_0))$ such that
     $$
    \theta_0 =   \angle_p( \vec u_0, \tilde q_{mid}) \le \frac{\pi}{2}. \tag2.17
     $$

     Let $\hat p = \phi(s_0)$, $\hat s_0 = d(p, \phi(s_0)) = |\vec u_0|$ and let $\sigma_{p\hat p}: [0, \hat s_0]
     \to M$ be a length-minimizing geodesic from $p$ to $p_0$ such that
     $$
\sigma_{p\hat p}'(0) = \frac{\vec u_0 }{|\vec u_0|}.
     $$

     We now apply inequality (2.8) to conclude that
     $$
   f( \phi(s_0) ) = f(\sigma_{p\hat p}(\hat s_0)) < f(p) - \hat s_0 \cos \theta_0 \le f(p).
     $$

     This completes the proof of the fact that $\Omega_{-c}$ is strictly convex.
   \qed
\enddemo

We emphasize that the equidistance evolution plays an important
role in Corollary 2.4. The conclusion of Corollary 2.4 fails if we replace
the distance function $f(x)=d(x,\partial\Omega_0)$ by an arbitrary
concave function $h:\Omega_0 \to [0,T_0]$.

For example, let $\Omega_0 =\{ (x,y)\in \Bbb{R}^2 :
x^2 +y^2 \leq 1 \}$ and $\Omega_{T_0} =\{ (x,0) : |x| \leq
\frac12 \}$. We can construct a convex function $\psi: \Omega_0
\to [-T_0,0]$ such that $\psi^{-1}(0) =\Omega_0$ and
$\psi^{-1} (-T_0) =\Omega_{-T_0}$. Choose $h(x)=-\psi(x)$. We notice
that $h^{-1} ([T,T_0]) =\Omega_{-T}$ is convex but not necessarily strictly
convex. For instance, $\Omega_{-T_0} =h^{-1} (T_0)$ is {\it not\/}
strictly convex and $\dim (\Omega_{-T_0})=1>0$. Thus, the conclusion
of Corollary 2.4 fails for non-distance function $h$.

We also remark that if the assumption of non-negative curvature is removed,
then the conclusion of Corollary 2.4 fails. Here is an example.

\bigskip
\noindent
{\bf Example 2.6.} Let $M^2 =\{ (x,y,z)| x^2 +y^2 -z^2 =1\}$ be a one-sheet hyperboloid in $\Bbb{R}^3$. It is clear that $M^2$ has negative
curvature. Let $\Omega_0 =\{ (x,y,z)\in M^2| z\leq 0\}$. It is clear that $\partial\Omega_0$ is a closed geodesic. Hence $\Omega_0$ is  a convex
subset of $M^2$. However, $\Omega_{-T} =\{ p\in \Omega_0 | d(p,\partial\Omega_0) \geq T\}$ has strictly concave boundary $\partial\Omega_{-T}$
for $T >0$. Thus $\Omega_{-T}$ is no longer a convex subset of $M^2$ for $T>0$.

\head $\S 3$ Further applications of sharp quadrangle comparisons
\endhead

In this section, we present two more applications of our new sharp
quadrangle comparisons.

First, we present a direct proof of Perelman's soul construction theorem.
His original proof used a contradiction argument.

\proclaim{Theorem 3.1}(Perelman [Per91,\S 6]) Let $M^n$ be a complete
Alexandrov space with $curv \geq k\geq 0$. Suppose that $\Omega_0$
is a convex sub-domain of $M^n$. Then $f(x)=d(x,\partial\Omega_0)$
is a concave function for $x\in\Omega_0$. Moreover, if $k>0$, then
$\Omega_c =f^{-1} ([c,+\infty))$ is strictly convex for $c>0$.
\endproclaim

\demo{Proof} We first consider the case of $k=0$, so $curv \geq 0$.
For each $p\in\Omega_0$ with $d(p,\partial\Omega_0)>0$ and for each
length-minimizing geodesic
$$
\sigma_{p\hat p} :[0,s] \to \Omega_0
$$
of unit speed, we would like to show that
$$
t\mapsto f(\sigma_{p\hat p}(t)) =\varphi_\sigma (t)
$$
is a concave function at $t=0$.

The derivative of $\varphi_\sigma$ is related to the angle $\theta$
between $\sigma'_{p\hat p} (0)$ and Min$(p,\partial\Omega_0)$,
where
$$
\align
\roman{Min} (p,\partial\Omega_0) &=\{ \sigma'_{pq}(0) |
\sigma_{pq}: [0,T] \to \Omega_0\
\text{is length-minimizing and of unit speed}  \\
&\text{from}\ p\ \text{to}\ q\in\partial\Omega_0 \ \text{with}\
d(p,q)=d(p,\partial\Omega_0)\}.
\endalign
$$
Let
$$
\theta =\min\{ \angle_p (\sigma'_{p\hat p}(0),\vec{w}) |
\vec{w}\in \roman{Min}(p,\partial\Omega_0) \}.
$$
We can show that
$$
\frac{d^+ \varphi_{\sigma}}{dt} (0) =-\cos\theta .
\tag 3.1
$$
By the proof of Corollary 2.4, we see that
$$
d(p,\partial\Omega_0) \leq d(p,q) -s\cos\theta,
\tag 3.2
$$
where $s=d(p,\hat p)$. It follows that
$$
f(\sigma_{p\hat p} (s)) \leq f(\sigma_{p\hat p} (0)) +
\frac{d^+ [f\circ \sigma_{p\hat p}]}{ds} (0) \, s.
\tag 3.3
$$
Therefore, $f(x)=d(x,\partial\Omega_0)$ is a concave function
for the case of $curv \geq 0$.
When $curv \geq k>0$, the proof of Corollary 2.4 implies that
$$
d(\hat p,\partial\Omega_0) <d(p,q)-s\cos\theta
\tag 3.4
$$
if $s>0$ and $0<\theta <\pi$. Equivalently, we have
$$
f(\sigma_{p\hat p} (s) )< f(\sigma_{p\hat p} (0) )+
\frac{d^+ [f\circ \sigma_{p\hat p}]}{ds} (0) \, s,
\tag 3.5
$$
when $s>0$ and $0<\theta <\pi$.
Notice that if $\sigma'_{p\hat p} (0)$ is tangent to
$\partial \Omega_c =f^{-1} (c)$ with $c=f(p)$, then $\theta =
\frac{\pi}{2}$. It follows from (3.5) that $\Omega_c$ is strictly convex for
any $c>0$.  \qed
\enddemo

Our second application is to improve Petrunin's second variational formula
for length in Alexandrov spaces of $curv \geq 0$.

Let us begin with convex curves in a $2$-dimensional Alexandrov space
with $curv \geq 0$.

\proclaim{Proposition 3.2} Let $M^2$ be a complete Alexandrov surface
with $curv  \geq 0$. Suppose that $\Omega_0$ is a convex sub-domain
of $M^2$ and $\Omega_c =\{ x\in\Omega_0 | d(x,\partial\Omega_0)
\geq c\}$. Then the length function $t\mapsto L(\partial\Omega_t)$
is a non-increasing function for $t\geq 0$.
\endproclaim

\demo{Proof} We have shown that $f(x)=d(x,\partial\Omega_0)$ is
a concave function. There is a Sharafutdinov semi-flow for gradient of $f$,
see [KPT07]. It was observed by Petrunin that the gradient flow
$$
\frac{d^+ \sigma}{dt} =\left. \frac{\nabla f}{|\nabla f|^2}
\right|_{\sigma(t)} \tag 3.6
$$
is a distance non-increasing map. Such a gradient flow induces a
Sharafutdinov projection $\partial\Omega_c \to \partial \Omega_{c+t}$
for $c\geq 0$ and $t\geq 0$. It is known that Sharafutdinov projection is
distance non-increasing. Thus we conclude that
$$
L(\partial\Omega_{c+t}) \leq L (\partial \Omega_c) \tag3.7
$$
for $c\geq 0$ and $t\geq 0$.
\qed
\enddemo

We now would like to elaborate the idea above. In fact, we can refine
the analysis above point-wise as follows. Let
$$
\pi_{c,c+t}:\partial \Omega_c \to \Omega_{c+t}
$$
be the Sharafutdinov distance non-increasing projection, where
$\dim [\Omega_c] =\dim (M^2)=2$. For each  curve
$\sigma_{c+t} \subset \partial\Omega_{c+t}$, we let
$$
\sigma_c =\pi^{-1}_{c,c+t} (\sigma_{c+t}).
$$

\proclaim{Corollary 3.4} Let $M^2$, $\Omega_0$, $\Omega_c$,
$\pi_{c,c+t}$, $\sigma_{c+t}$ and $\sigma_c$ be as above. Then
$$
L(\sigma_{c+t}) \leq L(\sigma_c)
\tag 3.8
$$
for any $c\geq 0$ and $t\geq 0$.
\endproclaim

We now consider the case when the initial curve $\sigma_0$ is a
geodesic segment. Let us give a short proof of Petrunin's 2nd
variational formula for this special case.

\proclaim{Theorem 3.5}(Improved Petrunin's formula for
$2$-dimensional case)
Let $M^2$ be a complete Alexandrov surface with non-negative curvature and $\Omega_0$
be a compact convex domain. Suppose that
$\hat\sigma: [-\varepsilon, L_0 +\varepsilon] \to M^2$ is a
length-minimizing geodesic of unit speed with
$\hat\sigma( [-\varepsilon, L_0 +\varepsilon] )\subset \partial
\Omega_0$,  $f(x)=d (x,\partial\Omega_0)$ and
$\varphi_q: [0,\delta] \to M^n$ is the gradient semi-flow
pointing inside $\Omega_0$ with
$$
\cases
\frac{d^+ \varphi_q}{dt} =\left. \frac{\nabla f}{|\nabla f|^2}
\right|_{\varphi_q (t)}, & \\
\varphi_q (0) =q\in\partial\Omega_0,  &
\endcases
$$
and suppose that $\sigma_t: [0, L_0] \to M^n$ is given by
$\sigma_t (s)=\varphi_{\hat\sigma (s)} (t)$. Then
$$
\frac{d^2 [L(\sigma_t)]}{dt^2} (0) \leq 0,
\tag 3.10
$$
where $L(\sigma_t)$ denotes the length of $\sigma_t$.
\endproclaim

\demo{Proof} Since $\sigma_0(s)=\hat\sigma (s)$ and
$\hat\sigma: [0,\ell] \to M^n$ is a length-minimizing geodesic,
we have
$$
\frac{d [L(\sigma_t)]}{dt} (0) =0.
\tag 3.11
$$
In addition, we have shown above that $t\mapsto L(\sigma_t)$
is a non-increasing function of $t$. Thus
$$
\frac{d [L(\sigma_t)]}{dt} (t) \leq 0.
\tag 3.12
$$
Therefore, we have
$$
\frac{d^2 [L(\sigma_t)]}{dt^2} (0) =\lim\limits_{t\to 0}
\frac{\frac{d^+ [L(\sigma_t)]}{dt}-\frac{d^+ [L(\sigma_0)]}{dt}}{t}
\leq 0.
$$
This completes the proof. \qed
\enddemo

The 2nd variational formula for lengths in higher dimensional Alexandrov
spaces will be discussed elsewhere. For the earlier work in this direction,
see Petrunin's paper [Petr98].

\head $\S 4$ Two open problems in Alexandrov's geometry and
possible approaches
\endhead

In this section, we discuss two open problems in Alexandrov's geometry
along with possible approaches. The first one is about the curvature
bound of the boundary of a convex domain in an Alexandrov space.
The second one is related to the geodesic semi-flow on Alexandrov spaces.

The following is a well-known problem in Alexandrov's geometry.

\proclaim{Open Problem 4.1} Let $M^n$ be a complete Alexandrov space with
$curv \geq k$ and $\Omega$ be a convex domain of $M^n$ with
$\dim (\Omega) =\dim (M^n)$. Prove that the boundary $\partial \Omega$
of $\Omega$ has $curv \geq k$ with respect to the intrinsic metric of
$\partial \Omega$.
\endproclaim

Earlier work for non-smooth convex domains in a smooth Riemannian manifold $M^n$
with $curv \geq k$ was carried out by Buyalo [Bu79]. Recently,
Alexander-Kapovitch-Petrunin ([AKP07]) showed that
$(\partial\Omega, d_{\partial\Omega})$ has $curv \geq k$ globally, i.e.
the conclusion of Toponogov comparison theorem holds for ``large" geodesic
triangles in $\partial\Omega$ as well.

 In Example 1.2 of an important paper [PP94], Perelman and Petrunin implicitly stated the following interesting result:
 {\it ``Let $\Omega_0$ be a convex
domain in a complete Alexandrov space $M^n$ with $curv \geq k$, and
$\overline{\Omega}_0$ be its closure in $M^n$. Suppose that
$\sigma:[0,\ell] \to\partial \Omega_0$ is a length-minimizing geodesic
segment of unit speed with respect to the intrinsic metric $d_{\partial
\Omega_0}$. The curve $\sigma$ is necessarily a quasi-geodesic segment
in $\overline{\Omega}_0$ (or in its doubling
$[\overline{\Omega}_0 \cup_{\partial
\Omega_0} \overline{\Omega}_0]$).}

One might take a different approach to the above Open Problem as follows.

\proclaim{Modified Problem 4.1.A} Prove that Toponogov comparison theorem holds for intrinsic
geodesic triangles in $\partial\Omega_0$.

    In particular, if $M^n$ has $curv \ge 0$ and if $\Delta \subset \partial\Omega_0$ is a triangle whose sides are
    length-minimizing segments with respect to the intrinsic metric of $\partial\Omega_0$, then is it
    true that the total (intrinsic) inner angles of $\Delta $ greater than or equal to $\pi$?
\endproclaim

In order to carry out our proposed approach above, we might want to
use an alternative definition of quasi-geodesics by introducing
generalized 2nd fundamental form for any subsets in an Alexandrov space
$M^n$.

Following Perelman-Petrunin's notion [Petr98], we say $\vec{v}=\log_p q$ if
there is a shortest path from $p$ to $q$ in $M^n$ which tangent to $\vec{v}$
and has length $|\vec{v}|=d_{M^n} (p,q)$. It is known that  $\log_p:
M^n \to T_p^- (M^n)$ is a distance non-decreasing map by ``global"
Toponogov comparison theorem.

For $q\not= p$ and $\vec{v} \in T_p^- (M^n)$ with $|\vec{v}|=1$, we let
$$
\angle_p (\vec{v},q) =\min\{ \angle_p (\vec{v},\vec{w}) |
\vec{w} \in \log_p (q) \}.
$$
For a $C^2$-smooth submanifolds $A$ in a smooth Riemannian manifold
$M^n$, the second fundamental form can be reviewed as follows.
If $\vec{v} \perp T_p (A)$ with $|\vec{v}|=1$ and if $\sigma:(-\varepsilon,
\varepsilon) \to A$ is a smooth curve of unit speed with $\sigma(0)=p$, then
$$
\align
\roman{II}_A^{\vec{v}} (\sigma'(0),\sigma'(0)) &=-\langle
\nabla_{\sigma'} \sigma' , \vec{v} \rangle |_{t=0} \\
&=\lim\limits_{\varepsilon\to 0} -\frac{1}{\varepsilon}
[ \langle \vec{v}, \exp_p^{-1} \sigma(\varepsilon) \rangle
+\langle \vec{v}, \exp_p^{-1} \sigma (-\varepsilon) \rangle ]. \tag 4.1
\endalign
$$

For example, let  $M^2=\Bbb{R}^2$ and
$A=S^1 =\{ (x,y)\in\Bbb{R}^2 | x^2 +y^2=1 \}$. We consider
$\sigma(t) =(\cos t,\sin t)$, $\vec{v}=(1,0)$. By the above formula, we have
$$
\roman{II}_A^{\vec{v}} (\sigma'(0),\sigma'(0)) =1>0.
$$
Inspired by discussion above and (4.1), we consider the generalized
2nd fundamental form for $A\subset M^n$ in the barrier sense.

\proclaim{Definition 4.2} Let $M^n$ be a complete Alexandrov space,
$p\in A\subset M^n$ and $\vec{v}\in T_p^- (M^n)$ with $|\vec{v}| =1$.
We let
$$
\theta_{p,\vec{v}} (\varepsilon) =\inf\{ \angle_p
(\vec{v},q) | q\in[A-B_\varepsilon (p)] \}.
$$
If
$$
\lim\limits_{\varepsilon\to 0^+} -\frac{\cos
\theta_{p,\vec{v}}}{\varepsilon} \leq 0
\tag 4.2
$$
holds, then we say that the subset $A$ is concave relative to $\vec{v}$
at $p$.
 \endproclaim

Using our new definition above, we are led to study the following problem.

\proclaim{Sub-Problem 4.1.B} (1) Let $\sigma: [0,\ell] \to M^n$
be a quasi-geodesic segment of unit speed in $M^n$ and
$A=\sigma([0,\ell])$. Prove that the quasi-geodesic $A=\sigma([0,\ell])$
is concave relative to any unit direction
$\vec{v}\in T_{\sigma(t)}^- (M^n)$.

(2) Let $\sigma: [0,\ell] \to M^n$ be a length-minimizing geodesic segment
of unit speed, $p=\sigma(t_0)$ with $0<t_0 <\ell$ and $\vec{v} \perp
\sigma'(t_0)$. Prove that
$$
\lim\limits_{\varepsilon\to 0^+}
\frac{\theta_{p,\vec{v}}(\varepsilon)}{\varepsilon} =0.
\tag 4.2
$$
\endproclaim

Our next question is related to the quasi-geodesic semi-flow on
Alexandrov spaces. In an earlier work of Perelman-Petrunin, quasi-geodesic
segments were extended in a ``non-unique" way. This was due the
fact that the choice of polar vectors are not unique.

\proclaim{Definition 4.3} Let $M^n$ be a complete Alexandrov space with
$curv \geq k$. Two unit tangent vectors $\{ \vec{v}_1, \vec{v}_2 \}
\subset T_p^- (M^n)$ are said to be polar if
$$
\cos \angle_p (\vec{v}_1,\vec{w}) +\cos \angle_p (\vec{v}_2,\vec{w})
\geq 0 \tag 4.4
$$
for any $\vec{w}\in T_p^- (M^n)$ with $|\vec{w}|=1$.
\endproclaim

Recall that
$$
\cos\theta_1 +\cos\theta_2 =2\cos \frac{\theta_1 +\theta_2}{2}
\cos \frac{\theta_1 -\theta_2}{2} .
\tag 4.5
$$
Because the unit tangent cone $\Sigma_p (M^n)$ has $curv
\geq 1$, its diameter is $\leq \pi$. The polar condition (4.4) is equivalent to
$$
\angle_p (\vec{v}_1,\vec{w}) + \angle_p (\vec{v}_2,\vec{w})
\leq \pi
\tag 4.6
$$
for all $\vec{w}\in \Sigma_p^- (M^n)$,

When the diameter of $\Sigma_p (M^n)$ is $\leq \frac{\pi}{2}$,
for any given $\vec{v}_1\in \Sigma_p(M^n)$, there are infinitely many
$\vec{v}_2 \in \Sigma_p(M^n)$ polar to $\vec{v}_1$.

However, if the radius of $\vec{v}_1$ in $\Sigma_p (M^n)$ satisfies
$$
\roman{rad} (\vec{v}_1) =\max \{ \angle_p
(\vec{v}_1,\vec{w}) | \vec{w}\in \Sigma_p (M^n) \} >
\frac{\pi}{2},\tag 4.7
$$
then there is a unique canonical choice of $\vec{v}_2$ which is polar
to $\vec{v}_1$ with $\angle_p (\vec{v}_1,\vec{v}_2) =
\roman{rad} (\vec{v}_1)$.

\proclaim{Proposition 4.5} Let $M^n$ be a complete Alexandrov space
with $curv \geq k$, $\vec{v}_1\in \Sigma_p(M^n)$ has
$\roman{rad} (\vec{v}_1) >\frac{\pi}{2}$. Then
$$
\Omega_{\frac{\pi}{2} +\varepsilon} =\{ \vec{w}\in T_p (M^n) |
\angle_p (\vec{v}_1,\vec{w}) \geq \frac{\pi}{2} +\varepsilon \}
$$
has strictly convex boundary $\partial  \Omega_{\frac{\pi}{2}
+\varepsilon}$ in $\Sigma_p(M^n)$ for any
$0< \varepsilon < [\roman{rad} (\vec{v}_1)-\frac{\pi}{2} ]$. Consequently,
there is a unique unit vector $\vec{v}_2 \in \Sigma_p(M^n)$ with
$\angle_p (\vec{v}_1, \vec{v}_2) =\roman{rad} (\vec{v}_1)$.
\endproclaim

\demo{Proof} In fact, $\vec{v}_2$ is the unique soul point of
$\Omega_{\frac{\pi}{2}+\varepsilon}$, where we used the fact that
$$
\roman{Hess} (f) \leq \cot (f) +df\otimes df
\tag 4.8$$
in the barrier sense, where $f(\vec{w})= d_{\Sigma} (\vec{v}_1,\vec{w})$.
Thus $\Omega_{\frac{\pi}{2}} =f^{-1} ([\frac{\pi}{2},\pi])$ is convex and
$\Omega_{\frac{\pi}{2}+\varepsilon} =f^{-1}([\frac{\pi}{2}+\varepsilon,
\pi ])$ is strongly convex for $\varepsilon >0$. \qed
\enddemo

In order to see how polar vectors are related to quasi-geodesics, we recall
a result of Perelman-Petrunin.

\proclaim{Proposition 4.6}([PPe 94])
Let $\sigma_1: [0,\ell_1] \to M^n$ and $\sigma_2: [\ell_1, \ell_2] \to M^n$
be two quasi-geodesic segments of unit speed with $\sigma_1 (\ell_1)
=\sigma_2 (\ell_1)$. Suppose that $\vec{v}_1=\frac{d^- \sigma_1}{dt}
(\ell_1)$ is polar to $\vec{v}_2=\frac{d^+ \sigma_2}{dt} (\ell_1)$, where
$$
\frac{d^{\pm} \sigma}{dt} (t) =\lim\limits_{\varepsilon \to 0^+}
\frac{\log_{\sigma(t)} \sigma(t\pm\varepsilon)}{\varepsilon} .
\tag 4.9
$$
Then $\sigma_1 \cup \sigma_2$ forms an extended quasi-geodesic segment,
where
$$
(\sigma_1 \cup \sigma_2)(t) =
\cases
\sigma_1 (t), & \roman{if}\ 0\leq t\leq \ell_1; \\
\sigma_2 (t), & \roman{if}\ \ell_1\leq t\leq \ell_2. \\
\endcases
$$
\endproclaim

Inspired by Propositions 4.4-4.6, we introduce the notion of canonical
quasi-geodesics.

\proclaim{Definition 4.7} Let $\sigma:[0,\ell] \to M^n$ be
a quasi-geodesic segment of unit speed. Suppose that
(1) $\roman{rad}\left( \frac{d^- \sigma}{dt} (t) \right)
>\frac{\pi}{2}$ and (2)
$\angle_{\sigma(t)} \left(  \frac{d^- \sigma}{dt} (t),
\frac{d^+ \sigma}{dt} (t)   \right) =
\roman{rad}\left( \frac{d^- \sigma}{dt} (t)\right)$ for $t\in (0,\ell)$.
Then $\sigma:[0,\ell] \to M^n$ is called a canonical quasi-geodesic.
\endproclaim

\proclaim{Definition 4.8} Let $M^n$ be a complete Alexandrov space with
$curv \geq k$. Let
$$
W^{reg} (M^n) =\{ (p,\vec{w}) | p\in M^n, \vec{w}\in \Sigma_p (M^n),
\roman{rad} (\vec{w}) > \frac{\pi}{2} \}
$$
be the non-extremal portion of unit cones over $M^n$. The we call
$W^{reg} (M^n)$ the canonical regular part of unit cone set over $M^n$.
\endproclaim

We conclude our paper by the following refined version of a problem of
Perelman and Petrunin.

\proclaim{Open Problem 4.9} Let $M^n$ be a complete Alexandrov space
with  $curv \geq k$. Suppose that $M^n$ has no boundary and $\vec{v}
\in W^{reg} (M^n)$ is a unit vector with $\roman{rad} (\vec{v})
> \frac{\pi}{2}$ in $\Sigma_p (M^n)$. Prove that there exists at most one
canonical quasi-geodesic $\sigma:[0,\ell] \to M^n$ with
$\sigma'(0) =\vec{v}$.

In addition, suppose that $U\subset W^{reg} (M^n)$ is a compact measurable subset
in regular part and that for each $(p,\vec{v}) \in U$, there is a canonical
quasi-geodesic $\varphi_{p,\vec{v}} :[0,\delta] \to M^n$ with
$$
\cases
\frac{d^+ \varphi_{p,\vec{v}} }{dt} (0) =\vec{v} & \\
\varphi_{p,\vec{v}} (0) =p; &
\endcases
$$
and $\psi_t (U) =\{ \varphi_{p,\vec{v}} (t) | (p,\vec{v}) \in U \}$.
Is it true that
$$
\roman{Vol}(\psi_t (U)) \leq \roman{Vol} (U) \tag4.10
$$
for $t\geq 0$?
\endproclaim

The volume non-increasing property (4.10) might be related to the concavity
properties of quasi-geodesics.

\medskip
\noindent{\bf Acknowledgement:}  There are some overlap between the recent work of Alexander-Bishop [AB08] and
Corollary 2.4 of our paper. Professors Stephanie Alexander kindly showed the first author the portion of  her manuscript
with Professor Bishop [AB08]  in October 2007. The overlap part of their work was earlier than ours. In addition, their work can also be applicable to
spaces with curvature bounded from above. Finally, we thank Professor Richard Bishop and Stephanie Alexander for pointing out some
inaccurate and ambiguous statements in an earlier version of our paper.


\enddocument